\documentclass[10pt,twoside]{article}
\usepackage{Latex-document}
\usepackage{multicol}
\usepackage{graphicx}
\def\volumeno{I}
\title{\bf The Work of Vladimir Voevodsky\vskip 6mm}
\author{Christophe Soul\'e\thanks{CNRS and Institut des Hautes
Etudes Scientifiques, 35, route de Chartres, 91440 Bures-sur-Yvette,
France.}\vspace*{-0.5cm}}
\date{\vspace{-8mm}}

\begin{document}

\maketitle \thispagestyle{first} \setcounter{page}{99}

\vskip 12mm

Vladimir Voevodsky was born in 1966. He studied at  Moscow State
University and  Harvard university. He is now Professor at the Institute
for Advanced Study in Princeton.

Among his main achievements are the following: he defined and developed
motivic cohomology and the ${\mathbf A}^1$-homotopy theory of algebraic
varieties; he proved the Milnor conjectures on the $K$-theory of fields.

Let us state the first Milnor conjecture. Let $F$ be a field and $n$ a
positive integer. The {\it Milnor $K$-group} of $F$ is the abelian group
$K_n^M (F)$ defined by the following generators and relations. The
generators are sequences $\{ a_1 , \ldots , a_n \}$ of $n$ units $a_i
\in F^*$. The relations are
\begin{eqnarray}
&&\{ a_1 , \ldots , a_{k-1} , xy , a_{k+1} , \ldots , a_n \} \nonumber \\
&= &\{ a_1 , \ldots , a_{k-1} , x , a_{k+1} , \ldots , a_n \} + \{ a_1 ,
\ldots , a_{k-1} , y , a_{k+1} , \ldots , a_n \} \nonumber
\end{eqnarray}
for all $a_i , x , y \in F^*$, $1 \leq k \leq n$, and the {\it Steinberg
relation}
$$
\{ a_1 , \ldots , x , \ldots , 1-x , \ldots , a_n \} = 0
$$
for all $a_i \in F^*$ and $x \in F - \{ 0,1 \}$.

On the other hand, let $\overline F$ be an algebraic closure of $F$ and
$G = {\rm Gal} (\overline F / F)$ the absolute Galois group of $F$, with
its profinite topology. The {\it Galois cohomology} of $F$ with
${\mathbf Z} / 2$ coefficients is, by definition,
$$
H^n (F , {\mathbf Z} / 2) = H_{\rm continuous}^n (G , {\mathbf Z} / 2)
\, .
$$

\vskip 6mm

{\bf Theorem 1.} (Voevodsky 1996 \cite{5}) \it  Assume $1/2 \in F$ and
$n \geq 1$. The Galois symbol
$$
h_n : K_n^M (F) / 2 K_n^M (F) \rightarrow H^n (F , {\mathbf Z} / 2)
$$
is an \rm isomorphism.

\vskip 6mm

This was conjectured by Milnor in 1970 \cite{1}. When $n=2$, Theorem~1
was proved by Merkurjev in 1983. The case $n=3$ was then solved
independently by Merkurjev-Suslin and Rost.

There exists also a Galois symbol on $K_n^M (F) / p \, K_n^M (F)$ for
any prime $p$ invertible in $F$. When $n=2$ and $F$ is a number field,
Tate proved that it is an isomorphism. In 1983 Merkurjev and Suslin
proved that it is an isomorphism when $n=2$ and $F$ is any field. Both
Voevodsky and Rost have made a lot of progress towards proving that, for
any $F$, any $n > 0$ and any $p$ invertible in $F$, the Galois symbol is
an isomorphism.

The map $h_n$ in Theorem~1 is defined as follows. When $n=1$, we have
$K_1^M (F) = F^*$ and $H^1 (F , {\mathbf Z} / 2) = {\rm Hom} (G ,
{\mathbf Z} / 2)$. The map
$$
h_1 : F^* / (F^*)^2 \rightarrow {\rm Hom} (G , {\mathbf Z} / 2)
$$
maps $a \in F^*$ to the quadratic character $\chi_a$ defined by
$$
\chi_a (g) = g (\sqrt a ) / \sqrt a = \pm 1
$$
for any $g \in G$ and any square root $\sqrt a$ of $a$ in $\overline F$.
That $h_1$ is bijective is a special case of Kummer theory. When $n \geq
2$, we just need to define $h_n$ on the generators $\{ a_1 , \ldots ,
a_n \}$ of $K_n^M (F)$. It is given by a cup-product:
$$
h_n (\{ a_1 , \ldots , a_n \}) = \chi_{a_1} \cup \cdots \cup \chi_{a_n}
\, .
$$
The fact that $h_n$ is compatible with the Steinberg relation was first
noticed by Bass and Tate.

Theorem~1 says that $H^n (F , {\mathbf Z} / 2)$ has a very explicit
description. In particular, an immediate consequence of Theorem~1 and
the definition of $h_n$ is the following

\vskip 6mm

{\bf Corollary 1.} \it The graded ${\mathbf Z} / 2$-algebra
$\displaystyle \bigoplus_{n \geq 0} H^n (F , {\mathbf Z} / 2)$ is
spanned by elements of degree one. \rm

\vskip 6mm

This means that absolute Galois groups are very special groups. Indeed,
it is seldom seen that the cohomology of a group
 or a topological space is spanned in
degree one.

\vskip 6mm

{\bf Corollary 2.} (Bloch) \it Let $X$ be a complex algebraic variety
and

\noindent $\alpha \in H^n (X ({\mathbf C}) ,$ ${\mathbf Z})$ a class in
its singular cohomology. Assume that $2 \alpha = 0$. Then, there exists
a nonempty Zariski open subset $U \subset X$ such that the restriction
of $\alpha$ to $U$ vanishes. \rm

\vskip 6mm

If Theorem~1 was extended to $K_n^M (F) / p \, K_n^M (F)$ for all $n$
and $p$, Corollary~2 would say that any torsion class in the integral
singular cohomology of $X$ is supported on some hypersurface. (Hodge
seems to have believed that such a torsion class should be Poincar\'e
dual to an analytic cycle, but this is not always true.)

With Orlov and Vishik, Voevodsky proved a second conjecture of Milnor
relating the Witt group of quadratic forms over $F$ to its Milnor
$K$-theory \cite{3}.

A very serious difficulty that Voevodsky had to overcome to prove
Theorem~1 was that, when $n=2$, Merkurjev made use of the algebraic
$K$-theory of conics over $F$, but, when $n \geq 2$, one needed to study
special quadric hypersurfaces of dimension $2^{n-1} - 1$. And it is
quite hard to compute the algebraic $K$-theory of varieties of such a
high dimension. Although Rost had obtained crucial information about the
$K$-theory of these quadrics, this was
 not enough to
conclude the proof when $n > 3$. Instead of algebraic $K$-theory,
Voevodsky used {\it motivic cohomology}, which turned out to be more
computable.

Given an algebraic variety $X$ over $F$ and two integers $p,q \in
{\mathbf Z}$, Voevodsky defined an abelian group $H^{p,q} (X , {\mathbf
Z})$, called motivic cohomology. These groups are analogs of the
singular cohomology of CW-complexes. They satisfy a long list of
properties, which had been anticipated by Beilinson and Lichtenbaum. For
example, when $n$ is a positive integer and $X$ is smooth, the group
$$
H^{2n,n} (X , {\mathbf Z}) = {\rm CH}^n (X)
$$
is the Chow group of codimension $n$ algebraic cycles on $X$ modulo
linear equivalence. And when $X$ is a point we have
$$
H^{n,n} ({\rm point}) = K_n^M (F) \, .
$$
It is also possible to compute Quillen's algebraic K-theory from motivic
cohomology. Earlier constructions of motivic cohomology are due to Bloch
(at the end of the seventies) and, later, to Suslin. The way Suslin
modified Bloch's definition was crucial to Voevodsky's approach and, as
a matter of fact, several important papers on this topic were written
jointly by Suslin and Voevodsky \cite{4,7}. There exist also two very
different definitions of $H^{p,q} (X, {\mathbf Z})$, due to Levine and
Hanamura; according to the experts they lead to the same groups. But it
seems fair to say that Voevodsky's approach to motivic cohomology is the
most complete and satisfactory one.

A larger context in which Voevodsky developed motivic cohomology is the
${\mathbf A}^1$-homotopy of algebraic manifolds \cite{6}, which is a
theory of ``algebraic varieties up to deformations'', developed jointly
with Morel \cite{2}. Starting with the category of smooth manifolds
(over a fixed field $F$), they first embed this category into the
category of {\it Nisnevich sheaves}, by sending a given manifold to the
sheaf it represents. A Nisnevich sheaf is a sheaf of sets on the
category of smooth manifolds for the Nisnevich topology, a topology
which is finer (resp. coarser) than the Zariski (resp. \'etale)
topology. Then Morel and Voevodsky define a homotopy theory of Nisnevich
sheaves in much the same way the homotopy theory of CW-complexes is
defined. The parameter space of deformations is the affine line
${\mathbf A}^1$ instead of the real unit interval $[0,1]$. Note that, in
this theory there are {\it two circles} (corresponding to the two
degrees $p$ and $q$ for motivic cohomology)! The first circle is the
sheaf represented by the smooth manifold ${\mathbf A}^1 - \{ 0 \}$
(indeed, ${\mathbf C} - \{ 0 \}$ has the homotopy type of a circle). The
second circle is ${\mathbf A}^1 / \{ 0,1 \}$ ( note that ${\mathbf R} /
\{ 0,1 \}$ is a loop).
 The
latter is not represented by a smooth manifold. But, if we identify $0$
and $1$ in the sheaf of sets represented by ${\mathbf A}^1$ we get a
presheaf of sets, and ${\mathbf A}^1 / \{ 0,1 \}$ can be defined as the
sheaf attached to this presheaf. This example shows why it was useful to
embed the category of algebraic manifolds into a category of sheaves.

It is quite extraordinary that such a homotopy theory of algebraic
manifolds exists at all. In the fifties and sixties, interesting
invariants of differentiable manifolds were introduced using algebraic
topology. But very few mathematicians anticipated that these ``soft''
methods would ever be successful for algebraic manifolds. It seems now
that any notion in algebraic topology will find a partner in algebraic
geometry. This has long been the case with Quillen's algebraic
$K$-theory, which  is precisely analogous to topological $K$-theory. We
mentioned that motivic cohomology is an algebraic analog of singular
cohomology. Voevodsky also computed the algebraic analog of the Steenrod
algebra, i.e. cohomological operations on motivic cohomology (this
played a decisive role in the proof of Theorem~1). Morel and Voevodsky
developed the (stable) ${\mathbf A}^1$-homotopy theory of algebraic
manifolds. Voevodsky defined {\it algebraic cobordism} as homotopy
classes of maps from the suspension of an algebraic manifold to the
classifying space MGL. There is also a direct geometric definition of
algebraic cobordism, due to Levine and Morel (see Levine's talk in these
proceedings), which should compare well with Voevodsky's definition. And
the list is growing: Morava $K$-theories, stable homotopy groups of
spheres, etc$\ldots$

Vladimir Voevodsky is an amazing mathematician. He
 has demonstrated an exceptional
talent for creating new abstract theories, about which he proved highly
nontrivial theorems. He was able to use these theories to solve several
of the main long standing problems in algebraic $K$-theory. The field is
completely different after his work. He opened large new avenues and, to
use the same word as Laumon,
 he is leading us
closer to the world of {\it motives} that Grothendieck was dreaming
about in the sixties.

\newpage

\title{\vspace*{12mm} \centerline{\Large \bf Vladimir Voevodsky}\vskip 6mm}
\author{}%\centerline{I.H.\'{E}.S., Bures-sur-Yvette, France}\vskip 1mm}
\date{\centerline{June 4, 1966}}
%N\'{e} le 6 novembre 1966 \`{a} Antony (Hauts-de-Seine), France} \vskip 1mm
%\centerline{\large Nationalit\'{e} fran\c{c}aise}}

\maketitle

\vskip 12mm

\begin{tabular}{p{2.5cm}p{8cm}}
  1989 & B.S. in Mathematics, Moscow University \\
  1992 & Ph.D. in Mathematics, Harvard University \\
  1992--1993 & Institute for Advanced Study, Member \\
  1993--1996 & Harvard University, Junior Fellow of Harvard Society of
  Fellows \\
  1996--1997 & Harvard University, Visiting Scholar \\
  1996--1997 & Max-Planck Institute, Visiting Scholar \\
  1996--1999 & Northwestern University, Associate Professor \\
  1998--2001 & Institute for Advanced Study, Member \\
  2002 & Institute for Advanced Study, Professor
\end{tabular}

\vfill

\begin{center}
\includegraphics[scale=0.7]{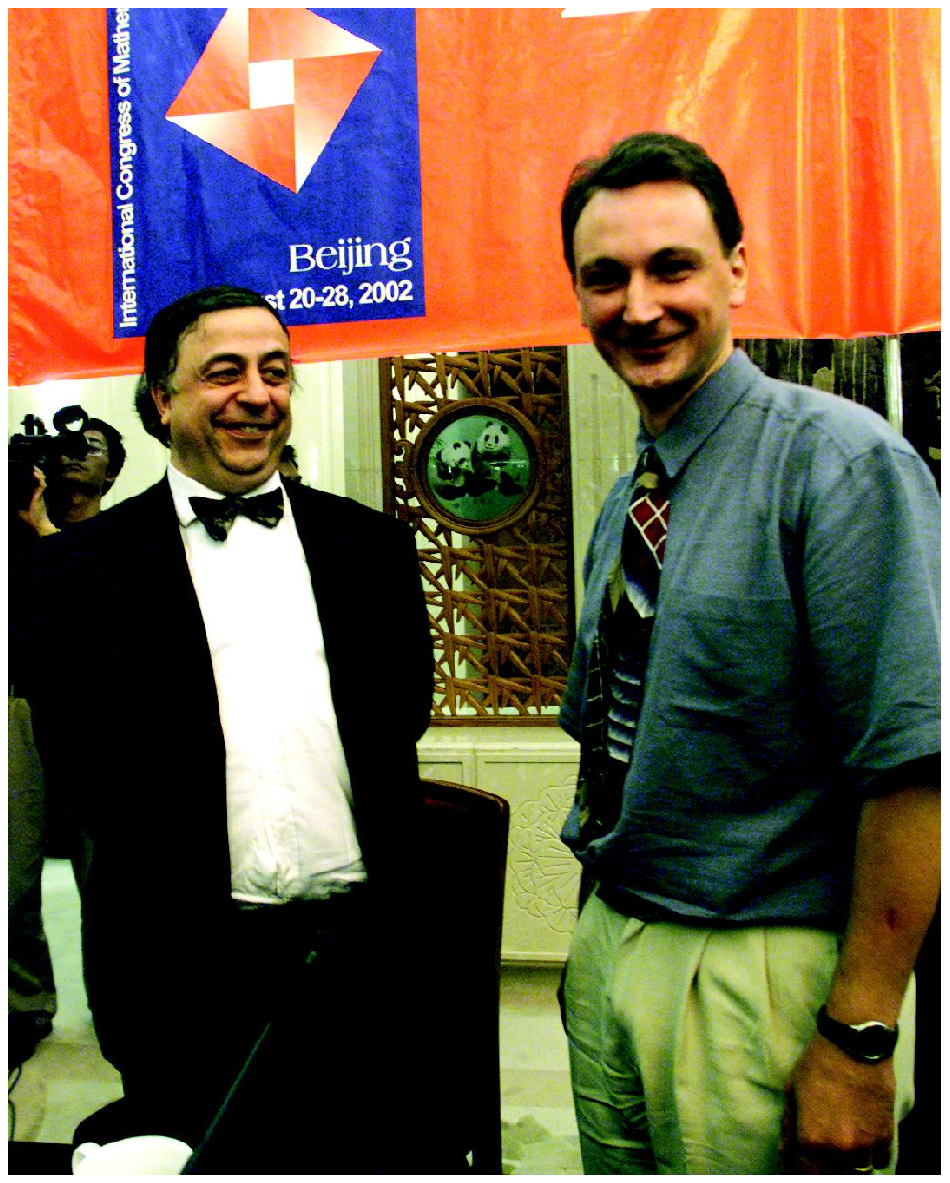}

\bigskip

V. Voevodsky (right) and C. Soul\'e
\end{center}

\label{lastpage}

\end{document}